\documentclass[10pt,draft
]{article}

\usepackage[latin1]{inputenc}
\usepackage{amsfonts}
\usepackage{amsmath}
\usepackage{amssymb}
\usepackage{amscd}
\usepackage{amsthm}
\usepackage{indentfirst}
\usepackage[hmargin=3.6cm
,vmargin=4.1cm
]{geometry}

\usepackage{epsfig}

\newtheorem{thm}{Theorem}[section]
\newtheorem{lemma}[thm]{Lemma}
\newtheorem{prop}[thm]{Proposition}
\newtheorem{coroll}[thm]{Corollary}

\theoremstyle{definition}

\newtheorem{defin}[thm]{Definition}
\newtheorem{rem}[thm]{Remark}
\newtheorem{exam}[thm]{Example}

\newtheorem*{acknow}{Acknowledgements}
\newtheorem*{prf}{Proof}

\newcommand{\R}{{\mathbb{R}}}

\newcommand{\T}{{\mathbb{T}}}
\newcommand{\Z}{{\mathbb{Z}}}

\newcommand{\C}{{\mathbb{C}}}

\newcommand{\cA}{{\mathcal{A}}}

\newcommand{\cG}{{\mathcal{G}}}

\newcommand{\cL}{{\mathcal{L}}}
\newcommand{\cM}{{\mathcal{M}}}

\newcommand{\cO}{{\mathcal{O}}}

\newcommand{\fc}{{:\ }}

\newcommand{\ol}{\overline}

\DeclareMathOperator{\IM}{Im}

\DeclareMathOperator{\Crit}{Crit}

\DeclareMathOperator{\im}{im}
\DeclareMathOperator{\id}{id}

\DeclareMathOperator{\osc}{osc}

\DeclareMathOperator{\Spec}{Spec}

\title{On the Lagrangian Hofer geometry in symplectically aspherical manifolds}
\author{Frol Zapolsky}

\begin{document}

\renewcommand{\labelenumi}{(\roman{enumi})}

\maketitle

\begin{abstract}
We use spectral invariants in Lagrangian Floer theory in order to show that there exist \emph{isometric} embeddings of normed linear spaces (finite or infinite dimensional, depending on the case) into the space of Hamiltonian deformations of certain Lagrangian submanifolds in tame symplectically aspherical manifolds. In addition to providing a new class of examples in which the Lagrangian Hofer metric can be computed explicitly, we refine and generalize some known results about it.
\end{abstract}

\section{Introduction and results}

Consider the following definition. Let $(W,\omega)$ be a symplectic manifold and let $L \subset W$ be a closed connected Lagrangian submanifold. Let $\cL(L)$ be the space of all Lagrangian submanifolds of $W$ Hamiltonian isotopic to $L$. For $L',L'' \in \cL$ we define
$$\rho(L',L'') = \inf\bigg\{\int_0^1\osc H_t\,dt\,\Big|\, \phi_H(L')=L''\bigg\}\,,$$
where $H_t$ is a compactly supported time-dependent Hamiltonian on $W$ and $\phi_H$ is the time-$1$ map of its flow. Chekanov \cite{Chekanov_Invariant_Finsler_metrics_Lagrangians} showed that in case $W$ is tame,\footnote{The tameness condition is needed for the compactness of spaces of holomorphic curves. \emph{Ibid.}, Chekanov gives an example of Lagrangian submanifold in a non-tame symplectic manifold for which $\rho$ is degenerate.} this quantity is a metric on $\cL(L)$.

Before stating our main result, we define the class of examples we deal with in this paper.
\begin{defin}
Let $(W,\omega)$ be a symplectic manifold. A closed connected Lagrangian submanifold $L\subset W$ is called relatively symplectically aspherical if $\omega|_{\pi_2(W,L)}=\mu|_{\pi_2(W,L)} = 0$, where $\mu$ is the Maslov index. A pair $(L,L')$ of closed connected Lagrangian submanifolds of $W$ is called weakly exact if for any smooth $u \fc S^1\times[0,1] \to W$ with $u(S^1\times 0) \subset L$, $u(S^1\times 1) \subset L'$ we have $\int u^*\omega = 0$.
\end{defin}
\begin{rem}Note that weak exactness of a pair of Lagrangians does not include a requirement on Maslov indices. Also note that if $L\subset W$ is relatively symplectically aspherical, then $W$ is necessarily symplectically aspherical in the sense of \cite{Schwarz_Action_sp_aspherical_mfds}, that is $\omega|_{\pi_2(W)}=c_1|_{\pi_2(W)}=0$.
\end{rem}

\noindent The main result of this paper is the following.
\begin{thm}\label{thm_main_result_1}
Assume that $(W,\omega)$ is tame and $L,L_0,\dots,L_k$ are relatively symplectically aspherical Lagrangian submanifolds, such that each pair $(L,L_j)$ is weakly exact, and such that $L$ intersects each $L_j$ transversely at a single point, and $L_j \cap L_{j'} = \varnothing$ for $j \neq j'$. Then there is an isometric embedding $(\R^k,\osc_0)\to(\cL(L),\rho)$.
\end{thm}
\noindent Here $\osc_0$ is a norm on $\R^k$ obtained as the restriction to $\R^k = 0\oplus \R^k \subset \R^{k+1}$ of the oscillation $\osc \fc \R^{k+1} \to [0,\infty)$ defined by $\osc (\tau_0,\dots,\tau_k) = \max_{j,j'}(\tau_j - \tau_{j'})$.

\begin{exam}
Examples of such manifolds are provided by the plumbing construction. Name\-ly, let $L,L_0,\dots,L_k$ be closed connected manifolds and consider their plumbing $W$ where the incidence graph is given by the tree with $L$ as the root and $L_0,\dots,L_k$ as the leaves.
\end{exam}

\begin{rem}
The novelty of this result is threefold. First, as mentioned in \cite{Usher_Hofer_metrics_bd_depth}, in case $W=T^*Q$, where $Q$ is closed and connected, there is an isometric embedding of $(C^\infty(Q)/\R,\osc)$ into $\cL(\cO_Q)$ where $\cO_Q \subset T^*Q$ is the zero section. However, the whole cotangent bundle is needed even in order to show that $\cL(\cO_Q)$ has infinite diameter. In contrast, in our result the manifold $W$ may be a small neighborhood of the union $L\cup\bigcup_jL_j$, in particular, it may have finite volume or symplectic capacities. Intuitively, in the cotangent bundle, the zero section can be ``moved'' in different directions inside the space $\cL(\cO_Q)$ along the various cotangent fibers. In our situation $L$ is ``moved'' in different direction along the $L_j$. Next, as remarked in \cite{Khanevsky_Hofer_metric_diameters}, sometimes one may pass to a covering space of the symplectic manifold and use the energy-capacity inequality to deduce results about the Hofer metric. In our situation $W$ may well be simply connected, so this method will not apply. Thirdly, in \cite{Usher_Hofer_metrics_bd_depth} Usher constructed quasi-isometric embeddings of normed linear spaces into $\cL$. Altough in certain cituations his embeddings are isometric, our theorem covers a new class of examples of such embeddings, in particular since $W$ may be non-compact.
\end{rem}

The next result covers some known cases, and also provides a new class of examples.
\begin{thm}\label{thm_main_result_2}
Let $L,L'$ be relatively symplectically aspherical Lagrangian submanifolds in $W$, such that the pair $(L,L')$ is weakly exact, and such that $L$ intersects $L'$ transversely at a single point. If in addition $L'$ admits a non-singular closed $1$-form\footnote{This is equivalent to $L'$ fibering over $S^1$.}, then there is an isometric embedding of $(C^\infty_c(0,1),\osc)$ into $(\cL(L),\rho)$.
\end{thm}
\begin{exam}
If $W$ is a surface, compact or not, and $L,L'$ are a pair of non-contractible curves intersecting at one point, then the theorem applies. This reproduces part of the result of Usher \cite[theorem 1.3]{Usher_Hofer_metrics_bd_depth}. As will be clear from the proof, when $W = \T^2$ and $L,L'$ is a couple of meridians intersecting at one point, there is actually an isometric embedding of $(C^\infty(S^1)/\R,\osc)$ into $\cL(L)$. More generally, if $W=\T^{2n}$ and $L,L'$ are two linear Lagrangian tori, the theorem applies, and as it is clear from the proof, there is an isometric embedding of $(C^\infty(\T^n)/\R,\osc)$ into $\cL(L)$.
\end{exam}
\begin{exam}
This theorem can be applied to plumbings, as follows. If $L,L'$ are two closed connected manifolds, we can form their single-intersection plumbing $W$, or more generally, we can take as $W$ any tame symplectic manifold which contains such a plumbing of $L,L'$, such that in $W$ both Lagrangians are relatively symplectically aspherical and such that $(L,L')$ is weakly exact.If $L'$ in addition admits a closed non-singular $1$-form, the theorem applies.
\end{exam}

\subsection{Preliminaries and notations}\label{section_preliminaries}

In the rest of the paper all manifolds are connected and all Lagrangian submanifolds are closed.

Throughout $(W,\omega)$ is a tame symplectic manifold of dimension $2n$. A Hamiltonian on $W$ is a function\footnote{When we concatenate Hamiltonians, we obtain functions in $C^\infty_c([0,2]\times W)$. The necessary adjustments are left to the reader.} $H \in C^\infty_c([0,1]\times W)$. If $W$ is closed, then, unless stated otherwise, all Hamiltonians are normalized, that is $\int_W H_t \omega^n = 0$ for all $t$. The Hamiltonian vector field of $H$ is defined via $\omega(X_{H_t},\cdot) = -dH_t$. The flow $\phi_H^t$ of $H$ is defined by $\phi_H^0 = \id$, $\frac{d\phi_H^t}{dt}=X_{H_t} \circ \phi_H^t$. We abbreviate $\phi_H = \phi_H^1$. The set of all time-$1$ flows of all the Hamiltonians is a group $\cG$, called the Hamiltonian group of $W$. When $L$ is a Lagrangian submanifold, the subgroup $\cG_L \subset \cG$ consists of all the diffeomorphisms fixing $L$ as a set.

When $K$ is a Hamiltonian, the reverse Hamiltonian $\ol K$ is defined via $\ol K(t,x) = -K(1-t,x)$. It generates the reverse isotopy $\phi_{\ol K}^t = \phi_K^{1-t}\phi_K^{-1}$. If $\gamma \fc [0,1] \to W$ is a path, the reverse path $\ol\gamma$ is defined by $\ol\gamma(t)=\gamma(1-t)$.

We will need to concatenate Hamiltonians and paths. If $z,z' \fc [0,1] \to Z$ are paths, where $Z$ is a space (in this paper $Z$ can be $W$, $\cG$, or $C^\infty_c(W)$) we let
$$(z\sharp z')(t)=\left\{\begin{array}{ll}z(t)\,, & t\leq 1\\ z'(t-1)\,, & t > 1\end{array}\right.\,.$$
This is a path defined over the interval $[0,2]$. It is continuous if $z(1)=z'(0)$ and smooth if all the derivatives of $z$ at $1$ coincide with those of $z'$ at $0$. A particularly important case is when we are given two Hamiltonians $H,K$ on $W$ and we need to concatenate them. There is a procedure, called smoothing, which allows us to do this for any pair of Hamiltonians, see for example \cite{Monzner_Vichery_Zapolsky_partial_qms_qss_cot_bundles}. Briefly, one chooses a function $f \fc [0,1] \to [0,1]$ which is smooth, monotone, and equals to $0$ near $0$ and to $1$ near $1$. Put $H^f(t,x)=f'(t) H(f(t),x)$. Then $H^f$ is also normalized, moreover it equals $0$ for $t$ near $0,1$. Therefore the concatenation $H^f\sharp K^g$, where $g$ is another smoothing function like this, is always well-defined and smooth. Moreover, dynamical invariants of $H$ are left intact, namely, the action spectrum, the set of periodic orbits or orbits with given boundary conditions, and the spectral invariants. More details can be found in \cite{Monzner_Vichery_Zapolsky_partial_qms_qss_cot_bundles}. Below, whenever we concatenate Hamiltonians, it is implicitly assumed that they have been previously smoothed.
                                                                                                                                                                                                                                                                                                                                                                                                                                                                 
We will also use various action functionals. Fix a Hamiltonian $H$. Let $\gamma$ be a smooth contractible loop in $W$ and let $u$ be a contracting disk, that is a map $u \fc D^2 \to W$ such that $\gamma = u|_{\partial D}$. Define
$$\cA_H(\gamma) = \int_0^1H_t(\gamma(t))\,dt - \int_{D^2} u^*\omega\,.$$
If $\omega|_{\pi_2(W)}=0$, then this is independent of $u$.

If $\gamma$ is a smooth path with endpoints on a Lagrangian $L$, we call it contractible for brevity, if $[\gamma]$ is the trivial element in $\pi_1(W,L)$. Let $u$ be a contraction of $\gamma$ into $L$, that is $u \fc D^2_+ \to M$, where $D^2_+ = \{z \in \C\,|\,|z|\leq 1\,,\IM z \geq 0\}$, with $u(\exp(\pi i t))=\gamma(t)$ and $u(t) \in L$ for $t\in[-1,1]$. Then put
$$\cA_H^L(\gamma) = \int_0^1H_t(\gamma(t))\,dt - \int_{D^2_+} u^*\omega\,.$$
When $\omega|_{\pi_2(W,L)}=0$, this is independent of $u$.

\subsection{Proofs of main results}

After the existence of spectral invariants in Lagrangian Floer homology has been established in the next section, the proof of main results is elementary, therefore we provide it here.

Before passing to the proof, we describe the construction of the so-called action homomorphism. It is instrumental for the methods of the present paper. Let $L \subset W$ be a relatively symplectically aspherical Lagrangian submanifold. If $\alpha \in \cG_L$ and $H$ is a Hamiltonian generating $\alpha$, we put $\gamma_q(t) = \phi_H^t(q)$ for $q\in L$. We have
\begin{prop}\label{prop_action_hom}
The Hamiltonian chords $\gamma_q$ are all contractible and the action $\cA_H^L(\gamma_q)$ only depends on $\alpha$. The map $\cA^L \fc \cG_L \to \R$ thus defined is a homomorphism. Moreover, if $H$ is a time-dependent Hamiltonian with $H_t|_L = c \in \R$ for all $t$, we have $\cA^L(\phi_H) = c$.
\end{prop}
\noindent The proof is given in \ref{section_Ham_loops_action_hom}. We call this $\cA^L$ the action homomorphism (associated to $L$).

Let us now formulate the properties of spectral invariants necessary for the proof of the main results.

\begin{thm}\label{thm_properties_sp_invts_intro}Let $L,L' \subset W$ be relatively symplectically aspherical Lagrangian submanifolds, which intersect transversely at a single point, and such that the pair $(L,L')$ is weakly exact. Then there is a function $\ell(\cdot:L,L') \fc \cG \to \R$ which satisfies:
\begin{enumerate}
\item $\int_0^1\min(F_t-G_t)\,dt \leq \ell(\phi_F:L,L') - \ell(\phi_G:L,L') \leq \int_0^1\max(F_t-G_t)\,dt$;
\item for $\alpha \in \cG_{L'}$ we have $\ell(\alpha\phi:L,L')=\ell(\phi:L,L') + \cA^{L'}(\alpha)$;
\item for $\beta \in \cG_L$ we have $\ell(\beta:L,L') = \cA^L(\beta)$; if $H$ is a Hamiltonian with $H_t|_L=c \in \R$, then $\ell(\phi_H:L,L') = \cA^L(\phi_H)=c$.
\end{enumerate}
\end{thm}
\noindent Theorem \ref{thm_properties_sp_invts_intro} follows from the more general result, theorem \ref{thm_sp_invts_summary}, to whose proof is dedicated most of section \ref{section_sp_invts_LFH}. Now we are ready to prove the main results. In order to compute Hofer distances between various Lagrangians, lower and upper bounds need to be established. Upper bounds are obtained via the oscillation of certain Hamiltonians and are elementary. The nontrivial part is to prove lower bounds. Since the technical idea behind the existence of such lower bounds is the same for both of the theorems, and is interesting on its own, we formulate it as a separate lemma.
\begin{lemma}\label{lemma_lower_bounds_Hofer_dist}
Let $Q,Q',Q''$ be relatively symplectically aspherical Lagrangian submanifolds of a tame symplectic manifold $(W,\omega)$, such that $Q,Q''$ are disjoint, $Q\cap Q'$ and $Q''\cap Q'$ are both transverse intersections which are single points, and such that both the pairs $(Q,Q')$, $(Q'',Q')$ are weakly exact. If $H$ is a Hamiltonian which satisfies $H|_Q=c$ and $H|_{Q''}=c''$, where $c,c'' \in \R$, then $\rho(\phi_H(Q'),Q') \geq c-c''$.
\end{lemma}
\begin{prf}
Abbreviate $\phi=\phi_H$. Theorem \ref{thm_properties_sp_invts_intro} implies that $\ell(\phi:Q,Q')=c$. Now let $\alpha \in \cG_{Q'}$ and let $G$ be a Hamiltonian generating $\alpha\phi$. It follows that
$$\int_0^1 \max G_t \,dt \geq \ell(\alpha\phi:Q,Q') = \ell(\phi:Q,Q') + \cA^{Q'}(\alpha) = c + \cA^{Q'}(\alpha)\,.$$
Analogously, $\ell(\phi:Q'',Q') = c''$, and also
$$-\int_0^1\min G_t\,dt \geq - \ell(\alpha\phi:Q'',Q') = -c'' - \cA^{Q'}(\alpha)\,.$$
The above two inequalities added together imply
$$\int_0^1\osc G_t\,dt \geq c-c''\,.$$
Taking infimum over $\alpha \in \cG_{Q'}$, we obtain
$$c-c'' \leq \rho(\phi(Q'),Q')\,. \qed$$ 
\end{prf}

\begin{prf}[of theorems \ref{thm_main_result_1} and \ref{thm_main_result_2}]

Let us first prove theorem \ref{thm_main_result_1}. Recall that we have the Lagrangians $L,L_0,\dots,L_k$. For $j=1,\dots,k$ let $H_j \in C^\infty_c(W,[0,1])$ be an autonomous Hamiltonian taking the value $1$ on $L_j$, and the value $0$ on all $L_{j'}$ for $j'=0,\dots,k$, $j' \neq j$, such that the supports of $H_j$ are all pairwise disjoint. For $\tau = (\tau_1,\dots,\tau_k)\in\R^k$ let $H_\tau = \sum_j\tau_j H_j$. Put $L_\tau = \phi_{H_\tau}(L)$. We claim that
$$\rho(L_\tau,L_{\tau'})=\osc_0(\tau-\tau')\,,$$
which is another way of saying that $(\R^k,\osc_0) \to (\cL(L),\rho),\,\tau \mapsto L_\tau$, is an isometric embedding. This will imply the assertion of the theorem.

Let us prove this. Abbreviate $\phi_\tau=\phi_{H_\tau}$ and note that
$$\rho(L_\tau,L_{\tau'}) = \rho(\phi_\tau(L),\phi_{\tau'}(L))=\rho(\phi_{\tau-\tau'}(L),L)\,.$$
Therefore it is enough to show that
$$\rho(L_\tau,L)=\osc_0(\tau)\,.$$
First,
$$\rho(L_\tau,L)=\rho(\phi_\tau(L),L)\leq\osc H_\tau\,.$$
It follows from the definition of $H_\tau$ that $\osc H_\tau = \osc_0(\tau)$. On the other hand, lemma \ref{lemma_lower_bounds_Hofer_dist} with $Q' = L$, $Q=L_j$, $Q'' = L_{j'}$, and $H = H_\tau$, implies
$$\rho(L_\tau,L) = \rho(\phi_\tau(L),L) \geq \tau_j - \tau_{j'}$$
for all $j,j' = 0,\dots,k$, where $\tau_0:=0$. Taking now the maximum over $j,j'$ results in
$$\osc_0(\tau) = \max_{j,j'=0,\dots,k}(\tau_j - \tau_j') \leq \rho(L_\tau,L) \leq \osc_0(\tau)\,,$$
which finishes the proof of theorem \ref{thm_main_result_1}.

Let us now prove theorem \ref{thm_main_result_2}. Recall that we have two Lagrangians $L,L'$. Let $\eta$ be a nonsingular closed $1$-form on $L'$. Identify a Weinstein neighborhood $U$ of $L'$ with a neighborhood of the zero section in $T^*L'$, such that $L \cap U = T_{q_0}^*L' \cap U$. Without loss of generality we assume that $U$ is the unit disk cotangent bundle of $L'$ with respect to some auxiliary Riemannian metric on $L'$. If necessary, scale $\eta$ so that its graph is contained in $U$. We let $L'_\tau$ be the Lagrangian in $W$ corresponding to the graph of $\tau\eta$ inside $U$, where $\tau\in[0,1]$. It follows that $L_\tau'$ is relatively symplectically aspherical and that the pair $(L'_\tau,L)$ is weakly exact, and moreover $L \cap L_\tau'$ is a single transverse intersection point.

Choose a function $H \in C^\infty_c(U)$ such that $H|_{L'_\tau}=\tau$ for $\tau \in [0,1]$. Extend $H$ by zero to $W$. Define a map $C^\infty_c(0,1) \to \cL(L)$ via $f \mapsto L_f:=\phi_f(L)$ where $\phi_f$ is the time-$1$ map of the Hamiltonian $H_f = f \circ H$. We claim that
$$\rho(L_f,L_{f'}) = \osc(f-f')\,,$$
which is what the assertion we want to establish says.
As in the first part, it is enough to show that $\rho(L_f,L) = \osc f$. On the one hand, we have
$$\rho(L_f,L) =\rho(\phi_f(L),L) \leq \osc H_f = \osc f\,.$$
On the other hand, if $\tau, \tau' \in (0,1)$, lemma \ref{lemma_lower_bounds_Hofer_dist} applied with $Q' = L$, $Q = L_\tau'$, $Q''=L_{\tau'}'$, and $H = H_f$, implies that
$$\rho(L_f,L) \geq f(\tau) - f(\tau')\,.$$
It follows that
$$\rho(L_f,L) \geq \max_{\tau,\tau' \in (0,1)} f(\tau) - f(\tau') = \osc f\,,$$
thereby completing the proof of the theorem. \qed
\end{prf}

\subsection{Discussion}

First we would like to remark that the existence of spectral invariants implies results on the Hofer geometry of the Hamiltonian group itself. For $\phi,\psi \in \cG$ we let
$$\rho(\phi,\psi)=\inf\left\{\int_0^1\osc H_t\,dt\,\Big|\,\phi\psi^{-1}=\phi_H\right\}\,.$$
This is a metric on $\cG$ (see for example \cite{Polterovich_geometry_of_Ham}). It follows from \cite{Leclercq_spectral_invariants_Lagr_FH} that whenever $L\subset W$ is a relatively aspherical Lagrangian, its Floer homology $HF(H:L)$ with Hamiltonian perturbations can be used to produce spectral invariants $\ell(A,\cdot:L)\fc \cG \to \R$ relative to classes $A \in HF(L) = H(L)$. These satisfy properties analogous to those of the relative invariants introduced below. See also \cite{Monzner_Vichery_Zapolsky_partial_qms_qss_cot_bundles}. In particular, the following can be easily deduced.
\begin{coroll}
Let $(W,\omega)$ be a tame symplectic manifold. (i) if $L_0,L_1,\dots,L_k \subset$ are relatively symplectically aspherical Lagrangians, all pairwise disjoint, then there is an isometric embedding $(\R^k,\osc_0) \to (\cG,\rho)$; (ii) if $L$ is a relatively symplectically aspherical Lagrangian which admits a non-singular closed $1$-form, then there is an isometric embedding $(C^\infty_c(0,1),\osc)\to (\cG,\rho)$. \qed
\end{coroll}

Secondly, we would like to point out that while the spectral invariants $\ell(\cdot:L,L')\fc \cG \to \R$ as in theorem \ref{thm_properties_sp_invts_intro} above are defined on the Hamiltonian group $\cG$, if $L''$ is another relatively symplectically aspherical Lagrangian, disjoint from $L$ and intersecting $L'$ transversely at one point, then the difference $\ell:=\ell(\cdot:L,L') - \ell(\cdot:L'',L')$ in fact descends to $\cL(L')$. The function thus obtained is in a sense (which can be made precise) a generalization of Viterbo's spectral invariants for Lagrangian submanifolds of cotangent bundles. This function $\ell\fc\cL(L') \to \R$ is Lipschitz in the Hofer metric: $\ell(K,K') \leq \rho(K,K')$ for $K,K' \in \cL(L')$. Moreover, lower bounds on this function are obtained as in lemma \ref{lemma_lower_bounds_Hofer_dist} via the action homomorphism. Namely, if $\beta \in \cG_L$ and $\beta'' \in \cG_{L''}$, then $\ell(\beta(L'),\beta''(L')) \geq \cA^L(\beta) - \cA^{L''}(\beta'')$. This is what makes it useful in computing Hofer's metric on $\cL(L')$.

\begin{acknow}
I would like to thank Kai Cieliebak and R\'emi Leclercq for helpful discussions, and Leonid Polterovich for useful suggestions. Part of this work was carried out during my stay at Ludwig-Maximilian-Universit\"at, Munich. I wish to thank this institution for a stimulating research atmosphere and hospitality.
\end{acknow}

\section{Spectral invariants in Lagrangian Floer theory}\label{section_sp_invts_LFH}

In this section we assume that $(W,\omega)$ is a tame symplectic manifold. We are given two transversely intersecting relatively symplectically aspherical Lagrangians $L,L'\subset W$, such that the pair $(L,L')$ is weakly exact. All homology is with coefficients in $\Z_2$ and the count of moduli spaces is modulo $2$.

\subsection{Definition and first properties}
We begin with a brief sketch of the construction of the Floer homology of the pair $(L,L')$ with Hamiltonian perturbations. General references are \cite{Floer_unregularized_grad_flow_symp_action} and \cite{Oh_FH_Lagr_intersections_hol_disks_I}.

We fix an intersection point $q_0 \in L \cap L'$ once and for all and consider the connected component $\Omega$ of the constant path at $q_0$ in the space of smooth paths $\{\gamma \fc [0,1]\to W\,|\, \gamma(0) \in L,\gamma(1) \in L'\}$. Since the pair $(L,L')$ is weakly exact, whenever $\gamma \in \Omega$ and $u$ is a homotopy from $q_0$ to $\gamma$, the integral
$$\int u^*\omega$$
is independent of $u$. Given a Hamiltonian $H$, we let $\cA_H^{L,L'} \fc \Omega\to\R$ be the action functional defined as
$$\cA_H^{L,L'}(\gamma) = \int_0^1 H_t(\gamma(t))\,dt - \int u^*\omega$$
for any homotopy $u$ from $q_0$ to $\gamma$. Its critical point set $\Crit \cA_H^{L,L'}$ consists precisely of those elements $\gamma \in \Omega$ which are Hamiltonian chords from $L$ to $L'$, that is which satisfy $\dot\gamma(t)=X_{H_t}(\gamma(t))$. The map $\Crit \cA_H^{L,L'} \to \phi_H(L)\cap L'$, $\gamma\mapsto \gamma(1)$ is injective. We let $CF(H:L,L')$ be the $\Z_2$-vector space\footnote{This vector space is ungraded since there is no requirement on Maslov indices.} spanned by $\Crit \cA_H^{L,L'}$. We call $H$ regular if $\phi_H(L)$ intersects $L'$ transversely. Regular Hamiltonians are generic. We choose $H$ regular, so that $\Crit \cA_H^{L,L'}$ is finite and $CF(H:L,L')$ is finite-dimensional.

Choose now a time-dependent compatible almost complex structure $J_t$ on $W$ for which $W$ is convex. This gives rise to an $L^2$-type metric on $\Omega$ and the negative gradient equation for $\cA_H^{L,L'}$ translates into Floer's PDE
$$\frac{\partial u}{\partial s} + J_t(u)\left(\frac{\partial u}{\partial t} - X_{H_t}(u)\right)=0$$
for $u \fc \R \times [0,1] \to W$ satisfying the boundary conditions $u (\R \times 0) \subset L$, $u(\R \times 1) \subset L'$. For $\gamma_\pm \in \Crit \cA_H^{L,L'}$ we let $\widehat{\cM}(\gamma_-,\gamma_+;H,J)$ be the space of solutions of this PDE subject to the asymptotic conditions $u(\pm\infty,\cdot)=\gamma_\pm$. For a generic choice of $J$ this is a finite-dimensional smooth manifold. Moreover, note that if $u \in \widehat{\cM}(\gamma_-,\gamma_+;H,J)$ then
$$E(u):=\int_{\R\times[0,1]}\left|\frac{\partial u}{\partial s}\right|^2\,ds\,dt = \cA_H^{L,L'}(\gamma_-) - \cA_H^{L,L'}(\gamma_+)\,.$$
Also note that $E(u) \geq 0$ with equality if and only if $u$ is a constant map. In case $\gamma_- \neq \gamma_+$ we let $\cM(\gamma_-,\gamma_+;H,J)$ be the quotient by the natural action of $\R$. We also put $\cM(\gamma_-,\gamma_-;H,J) = \varnothing$.

The Floer boundary operator on $CF(H:L,L')$ is given on generators by
$$\partial \gamma_- = \sum _{\gamma_+:\dim\cM(\gamma_-,\gamma_+;H,J)=0}\#\cM(\gamma_-,\gamma_+;H,J)\,\gamma_+\,,$$
where for a manifold $Y$ we let $\#Y$ be the modulo $2$ number of points in the zero-dimensional part of $Y$. It is a standard fact that $\partial^2 = 0$ and we let $HF(H:L,L')$ be the corresponding homology. We omit the almost complex structure from the notation since neither the Floer homology nor the spectral invariants which we will introduce shortly depend on it.

Since $L,L'$ intersect transversely, the zero Hamiltonian is regular. In this case $HF(0:L,L')$ is nothing but the (component at $q_0$ of the) usual Floer homology $HF(L,L')$ of the pair $(L,L')$. We record a particularly important special case:
\begin{lemma}
If $q_0$ is the only intersection point of $L,L'$ then $HF(0:L,L') = \Z_2\cdot q_0 \equiv \Z_2$. \qed
\end{lemma}
\noindent This is true since the boundary operator vanishes. Indeed, by definition $\cM(q_0,q_0;0,J) = \varnothing$.

For any two regular Hamiltonians $H,H'$ and any two regular almost complex structures $J,J'$ there is a canonical continuation isomorphism $HF(H:L,L') = HF(H':L,L') = HF(L,L')$. Since the boundary operator counts negative gradient lines of the action functional, it follows that the subspace $CF^a(H:L,L') \subset CF(H:L,L')$ generated by the critical points of $\cA_H$ of action $<a$ is a subcomplex. We let $i^a \fc HF^a(H:L,L') \to HF(H:L,L') = HF(L,L')$ be the morphism induced by the inclusion. The aforementioned continuation isomorphisms leave the maps $i^a$ intact if we only change $J$, therefore we omit it from the notation throughout. We are now in position to define the spectral invariants. For $A \in HF(L,L')-\{0\}$ we let
$$\ell(A,H:L,L') = \inf \{a \,|\,A \in \im i^a\}\,.$$
From the definition it follows that the spectral invariants are spectral, that is $\ell(A,H:L,L')$ belongs to the action spectrum $\Spec(H:L,L') = \cA_H^{L,L'}\big(\Crit \cA_H^{L,L'}\big)$.

We have defined spectral invariants for a regular Hamiltonian. The existence of the continuation isomorphisms implies the following:
$$\int_0^1 \min(H_t - H_t')\,dt \leq \ell(A,H:L,L') - \ell(A,H':L,L') \leq \int_0^1 \max(H_t - H_t')\,dt$$
for regular $H,H'$. These inequalities allow us to define $\ell(A,G:L,L')$ for an arbitrary smooth Hamiltonian\footnote{In fact, at the same price we can define them for arbitrary continuous Hamiltonians, however this will not be needed in this paper.} $G$, in the standard manner. Using techniques similar to those of \cite{Oh_Construction_sp_invts_Ham_paths_closed_symp_mfds}, we can show that these extended invariants are also spectral.

\subsection{Hamiltonian loops and the action homomorphism}\label{section_Ham_loops_action_hom}

In order to be able to use spectral invariants in estimates on the Lagrangian Hofer distance, we need to show that they are defined on the Hamiltonian group, that is, that $\ell(A,H:L,L')$ in fact only depends on the time-$1$ map $\phi_H$. The first step is to show that the action of any periodic orbit of a Hamiltonian generating a loop in $\cG$ is zero.

Since $L$ is relatively symplectically aspherical, $W$ is symplectically aspherical. Let $G$ be a Hamiltonian generating a loop\footnote{In agreement with our convention, since $G$ needs to be time-periodic in order for what follows to make sense, we smooth it as described in subsection \ref{section_preliminaries}.} in $\cG$. Consider first the case of $W$ closed. The existence of Floer homology (see for example \cite{Schwarz_Action_sp_aspherical_mfds}) implies that every $1$-periodic orbit of the form $t\mapsto \phi_G^t(z)$ is contractible. Moreover, in \cite{Schwarz_Action_sp_aspherical_mfds} it is proved that the action of every such orbit is zero (remember that $G$ is normalized). If $W$ is open, there is a simpler argument. First, all the periodic orbits of $G$ are homotopic, which implies that all of them are contractible, because orbits outside the support of $G$ are just constant. It also follows that the action of these constant orbits is zero, and since the action is independent of the beginning of the orbit, all the actions of the periodic orbits vanish. We thus have
\begin{lemma}Let $G$ be a Hamiltonian such that $\phi_G = \id$. Then any periodic orbit of the form $t \mapsto \phi_G^t(z)$, where $z\in W$, is contractible, its action is well-defined and equal to zero. \qed
\end{lemma}

We can now establish the existence of the action homomorphism.
\begin{prf}[or proposition \ref{prop_action_hom}]It follows from the existence of Floer homology for relatively symplectically aspherical Lagrangians (see for example \cite{Leclercq_spectral_invariants_Lagr_FH}) that at least one of the chords $\gamma_q$ is contractible. Since all of them are homotopic, it follows that all of them are contractible. This implies that the actions $\cA_H^L(\gamma_q)$ are well-defined and are all equal. It remains to show that this number is independent of $H$.

Put $\gamma_q^H(t) = \phi_H^t(q)$. Let $K$ be another Hamiltonian generating $\phi_H$, and let $\gamma_q^K(t) = \phi_K^t(q)$. We need to prove that $\cA_K^L(\gamma_q^K)=\cA_H^L(\gamma_q^H)$. Pick a contracting half-disk $u$ for $\gamma_q^H$. The concatenation $\gamma_q^H \sharp \ol{\gamma_q^K}$ is a periodic orbit of the flow of the concatenated Hamiltonian $H \sharp \ol K$. Since this Hamiltonian generates the identity map, the above lemma implies that there is a contracting disk $v$ for $\gamma_q^H \sharp \ol{\gamma_q^K}$ and that
$$\cA_{H\sharp \ol K}\big(\gamma_q^H \sharp \ol{\gamma_q^K}\big) = 0\,.$$
Writing out the action, we have
$$\cA_{H\sharp \ol K}\big(\gamma_q^H \sharp \ol{\gamma_q^K}\big) = \int_0^2\big(H\sharp\ol K\big)_t\big(\big(\gamma_q^H \sharp \ol{\gamma_q^K}\big)(t)\big)\,dt - \int v^*\omega\,.$$
Since $v$ provides a homotopy with fixed endpoints from $\gamma_q^K$ to $\gamma_q^H$, the concatenation $\ol v\sharp u$ is a contracting half-disk for $\gamma_q^K$ and therefore
\begin{align*}
\cA_H^L(\gamma_q^H) - \cA_K^L(\gamma_q^K) &= \int_0^1 H_t(\gamma_q^H(t))\,dt - \int u^*\omega - \int_0^1 K_t(\gamma_q^K(t))\,dt + \int (\ol v\sharp u)^*\omega\\
&=\int_0^2\big(H\sharp\ol K\big)_t\big(\big(\gamma_q^H \sharp \ol{\gamma_q^K}\big)(t)\big)\,dt - \int v^*\omega = \cA_{H\sharp \ol K}\big(\gamma_q^H \sharp \ol{\gamma_q^K}\big) = 0
\end{align*}
by the discussion above. It remains to show that if $H$ is such that $H|_L = c \in \R$ then $\cA^L(\phi_H) = c$. This follows from the fact that the flow of such a Hamiltonian preserves $L$ and thus the area part of the action of any Hamiltonian chord vanishes since it is completely contained in $L$, while the first part is seen to be equal to $c$. This finishes the proof of the proposition. \qed
\end{prf}

\subsection{The shift of spectral invariants and independence of isotopy}

In this subsection we establish the shift property of spectral invariants and use it to prove their independence of isotopy.
\begin{thm}\label{thm_shift_sp_invts}
Let $K,H$ be Hamiltonians such that $\phi_K \in \cG_{L'}$. Put $G_t = K_t + H_t \circ (\phi_K^t)^{-1}$. Then for any $A \in HF(L,L')$ we have
$$\ell(A,G:L,L') = \ell(A,H:L,L') + \cA^{L'}(\phi_K)\,.$$
\end{thm}
\begin{prf}
The proof is based on the so-called naturality isomorphism, see for example \cite{Barraud_Cornea_Lagrangian_intersections_Serre_sp_seq}, \cite{Leclercq_spectral_invariants_Lagr_FH}. We only sketch the construction; the interested reader is referred to these papers for details.

The map
$$\Crit \cA_H^{L,L'} \to \Crit \cA_G^{L,L'}\,,\quad \gamma \mapsto \gamma^K\,,$$
where $\gamma^K(t) = \phi_K^t(\gamma(t))$, is well-defined and is a bijection. Perhaps the non-trivial point here is that it maps elements of $\Omega$ back to $\Omega$. This can be seen as follows. The curve $\gamma^K$ is homotopic with fixed endpoints to the concatenation $\gamma\sharp\delta$ where $\delta(t) = \phi_K^t(\gamma(1))$. But $\delta$ is a Hamiltonian chord of $K$, whose flow maps $L'$ back to itself, and therefore it is contractible into $L'$ by proposition \ref{prop_action_hom}. This shows that $\gamma^K$ is homotopic to $q_0$ if $\gamma$ is.

Since the path $\{\phi_G^t\}_t$, which is generated by $G$, is homotopic with fixed endpoints to the concatenation $\{\phi_H^t\}_t\sharp\{\phi_K^t\phi_H\}_t$, which is generated by the concatenated Hamiltonian $H \sharp K$, it follows from a standard argument (see for example, \cite{Oh_Construction_sp_invts_Ham_paths_closed_symp_mfds}) that
$$\cA_{H\sharp K}^{L,L'}(\gamma\sharp\delta)=\cA_G^{L,L'}(\gamma^K)\,.$$
Since action is additive under concatenations and $\cA_K^{L'}(\delta)=\cA^{L'}(\phi_K)$ by proposition \ref{prop_action_hom}, we have
$$\cA_G^{L,L'}(\gamma^K) = \cA_H^{L,L'}(\gamma)+\cA^{L'}(\phi_K)\,,$$
that is, $\gamma \mapsto \gamma^K$ shifts the action by $\cA^{L'}(\phi_K)$.

Next we have the map
$$\widehat{\cM}(\gamma_-,\gamma_+;H,J) \to \widehat{\cM}(\gamma_-^K,\gamma_+^K;G,J^K)\,,\quad u \mapsto u^K\,,$$
where $u^K(s,t) = \phi_K^t(u(s,t))$ and $J^K_t = \phi^t_{K,*}\,J_t\,(\phi^t_{K,*})^{-1}$. This map is a diffeomorphism. Clearly it preserves the $\R$-action and so passes to a diffeomorphism $\cM(\gamma_-,\gamma_+;H,J) \to \cM(\gamma_-^K,\gamma_+^K;G,J^K)$. The net result is that $\gamma\mapsto\gamma^K$ extends to a chain isomorphism
$$CF^a(H:L,L') \to CF^{a+\cA^{L'}(\phi_K)}(G:L,L')$$
for any $a$, and the assertion of the theorem follows. \qed
\end{prf}

\begin{rem}\label{rem_action_spectrum}
Note that the proof of this theorem implies that the action spectrum $\Spec(H:L,L')$ only depends on $\phi_H$. This follows from the definition of the above bijection $\gamma \mapsto \gamma^K$ which preserves the actions if $K$ generates the identity map since then $\cA^{L'}(\phi_K) = 0$.
\end{rem}

We can now prove
\begin{prop}
The spectral invariant $\ell(A,H:L,L')$ only depends on $\phi_H \in \cG$.
\end{prop}
\begin{prf}Abbreviate $\ell(\cdot)\equiv \ell(A,\cdot:L,L')$. If $K$ is another Hamiltonian generating $\phi_H$, then the concatenation $H \sharp \ol K$ generates the identity map. Since $\cA^{L'}$ is a homomorphism to the reals, it maps the identity map to zero. It follows from theorem \ref{thm_shift_sp_invts} that
$$\ell(H\sharp\ol K\sharp K) = \ell(K) + \cA^{L'}(\phi_{H\sharp \ol K}) = \ell(K) +\cA^{L'}(\id) = \ell(K)\,.$$
It now suffices to note that
$$\ell(H) = \ell(H\sharp\ol K\sharp K)\,.$$
For a proof, see, for example \cite{Monzner_Vichery_Zapolsky_partial_qms_qss_cot_bundles}. \qed
\end{prf}

\subsection{Summary}

The following theorem summarizes the various properties of spectral invariants.
\begin{thm}\label{thm_sp_invts_summary}The spectral invariants descend to functions $\ell(A,\cdot:L,L') \fc \cG \to \R$; moreover, they satisfy:
\begin{enumerate}
\item $\ell(A,\phi:L,L') \in \Spec(\phi)$;\footnote{See remark \ref{rem_action_spectrum} for the definition of the action spectrum $\Spec(\phi)$.}
\item $\int_0^1\min(F_t-G_t)\,dt \leq \ell(A,\phi_F:L,L') - \ell(A,\phi_G:L,L') \leq \int_0^1\max(F_t-G_t)\,dt$;
\item for $\alpha \in \cG_{L'}$ we have $\ell(A,\alpha\phi:L,L')=\ell(A,\phi:L,L') + \cA^{L'}(\alpha)$;
\item in case $\#L \cap L' = 1$, let $A \in HF(L,L') = \Z_2$ be the generator; then for $\beta \in \cG_L$ we have $\ell(A,\beta:L,L') = \cA^L(\beta).$
\end{enumerate}
\end{thm}
\begin{prf}
The only statement remaining to be proved is (iv). Let $H$ generate $\beta$. By spectrality $\ell(A,\beta:L,L') = \cA_H^{L,L'}(\gamma)$ where $\gamma \in \Omega$ is a Hamiltonian chord from $L$ to $L'$. We have $\gamma(1) = \phi_H(\gamma(0))$. Since by assumption $\phi_H(L) = L$, we see that $\gamma(1) \in L$, that is $L$ is a Hamiltonian chord beginning and ending on $L$. By proposition \ref{prop_action_hom}, there is a contracting half-disk $u$ for $\gamma$ into $L$. It can be seen that this $u$ may be reparametrized to produce a homotopy from $q_0$ to $\gamma$. Therefore
$$\cA_H^{L,L'}(\gamma) = \cA_H^L(\gamma) = \cA^L(\beta)\,,$$
as claimed. The proof of the theorem is complete. \qed
\end{prf}
\noindent Theorem \ref{thm_properties_sp_invts_intro} follows from the one we just proved, therefore the main results are now proved as well.

\bibliographystyle{plain}
\bibliography{1}

\begin{thebibliography}{10}

\bibitem{Barraud_Cornea_Lagrangian_intersections_Serre_sp_seq}
Fran{\c c}ois Barraud and Octav Cornea.
\newblock Lagrangian intersections and the {S}erre spectral sequence.
\newblock {\em Annals of Mathematics}, 166:657--722, 2007.

\bibitem{Chekanov_Invariant_Finsler_metrics_Lagrangians}
Yuri Chekanov.
\newblock Invariant {F}insler metrics on the space of {L}agrangian embeddings.
\newblock {\em Mathematische Zeitschrift}, 234(3):605--619, 2000.

\bibitem{Floer_unregularized_grad_flow_symp_action}
Andreas Floer.
\newblock Unregularized gradient flow of the symplectic action.
\newblock {\em Communications on Pure and Applied Mathematics}, 41:775--813,
  1988.

\bibitem{Khanevsky_Hofer_metric_diameters}
Michael Khanevsky.
\newblock Hofer's metric on the space of diameters.
\newblock {\em Journal of Topology and Analysis}, 1(4):407--416, 2009.

\bibitem{Leclercq_spectral_invariants_Lagr_FH}
R\'emi Leclercq.
\newblock Spectral invariants in {L}agrangian {F}loer theory.
\newblock {\em Journal of Modern Dynamics}, 2(2):249--286, 2008.

\bibitem{Monzner_Vichery_Zapolsky_partial_qms_qss_cot_bundles}
Alexandra Monzner, Nicolas Vichery, and Frol Zapolsky.
\newblock Partial quasi-morphisms and quasi-states on cotangent bundles, and
  symplectic homogenization.
\newblock {\tt arXiv:1111.0287}.

\bibitem{Oh_FH_Lagr_intersections_hol_disks_I}
Yong-Geun Oh.
\newblock Floer cohomology of {L}agrangian intersections and pseudo-holomorphic
  disks, {I}.
\newblock {\em Communications on Pure and Applied Mathematics}, 46(7):949--993,
  1993.

\bibitem{Oh_Construction_sp_invts_Ham_paths_closed_symp_mfds}
Yong-Geun Oh.
\newblock Construction of spectral invariants of hamiltonian paths on closed
  symplectic manifolds.
\newblock In {\em The breadth of symplectic and {P}oisson geometry}, volume 232
  of {\em Progress in Mathematics}, pages 525--570. Birkh\"auser {B}oston,
  Boston, MA, 2005.

\bibitem{Polterovich_geometry_of_Ham}
Leonid Polterovich.
\newblock {\em The geometry of the group of symplectic diffeomorphisms}.
\newblock Lectures in mathematics {ETH} {Z}\"urich. Birkh\"auser {V}erlag,
  Basel, 2001.

\bibitem{Schwarz_Action_sp_aspherical_mfds}
Matthias Schwarz.
\newblock On the action spectrum for symplectically aspherical manifolds.
\newblock {\em Pacific Journal of Mathematics}, 193(2):419--461, 2000.

\bibitem{Usher_Hofer_metrics_bd_depth}
Michael Usher.
\newblock Hofer's metrics and boundary depth.
\newblock {\tt arXiv:1107.4599}.

\end{thebibliography}

\end{document}